\documentclass[12pt]{article} 
\usepackage{amsfonts} 
\usepackage{amsmath} 
\usepackage{amsthm} 
\usepackage{amssymb} 
\usepackage{mathrsfs} 
\usepackage{amstext} 
\usepackage{url} 

\newtheorem{Theorem}{Theorem}[section] 
 
\newtheorem{corollary}[Theorem]{Corollary}

\theoremstyle{remark}  
\newtheorem{remark}{Remark}[section]

\renewcommand{\Re}{\operatorname{Re}}  
\renewcommand{\Im}{\operatorname{Im}}

\renewcommand{\qed}{\rule{2.5mm}{2.5mm}}  
\newenvironment{Proof}{\noindent  
{\bf\underline{Proof:} }}  
{\hspace*{\fill}\qed\vskip1em}

\numberwithin{equation}{section}  
\input{tcilatex}  

\begin{document} 

\title{Learning Trigonometric Polynomials from Random Samples and 
Exponential Inequalities for Eigenvalues of\ Random Matrices} 
\author{Karlheinz Gr\"{o}chenig\footnotemark[1]  , Benedikt M.~P{\"o}tscher%
\footnotemark[2]  , Holger Rauhut\footnotemark[3]  } 
\date{January 29, 2010} 
\maketitle 

\begin{abstract} 
Motivated by problems arising in random sampling of trigonometric 
polynomials, we derive exponential inequalities for the operator norm of the 
difference between the sample second moment matrix $n^{-1}U^{\ast }U$ and 
its expectation where $U$ is a complex random $n\times D$ matrix with 
independent rows. These results immediately imply deviation inequalities for 
the largest (smallest) eigenvalues of the sample second moment matrix, which 
in turn lead to results on the condition number of the sample second moment 
matrix. We also show that trigonometric polynomials in several variables can 
be learned from $const\cdot D\ln D$ random samples. 

\end{abstract} 


\renewcommand{\thefootnote}{\fnsymbol{footnote}} 

\footnotetext[1]{%
Department of Mathematics, University of Vienna, Nordbergstrasse 15, A-1090 
Vienna. E-mail: karlheinz.groechenig@univie.ac.at.\newline 
K.~G.~was supported by the Marie-Curie Excellence Grant MEXT-CT
2004-517154.} 
\footnotetext[2]{%
Department of Statistics, University of Vienna, Universit\"atsstrasse 5, 
A-1010 Vienna. E-mail: benedikt.poetscher@univie.ac.at} 
\footnotetext[3]{%
Hausdorff Center for Mathematics, Endenicher Allee 60, D-53115 Bonn, 
Germany, E-mail:  rauhut@hcm.uni-bonn.de. \newline 
H.~R.~was supported by the Individual Marie-Curie Fellowship MEIF-CT 
2006-022811.}
\renewcommand{\thefootnote}{\arabic{footnote}} 

\section{Introduction} 

Let $U$ be a complex random $n\times D$ matrix with independent rows. The 
matrix of (non-centered) sample second moments is then given by $%
n^{-1}U^{\ast }U$. We provide exponential probability inequalities for the 
operator norm of the difference between the sample second moment matrix and 
its expectation. These results immediately imply deviation inequalities for 
the largest (smallest) eigenvalues of the sample second moment matrix. As a 
consequence we obtain probability inequalities for the condition number of 
the sample second moment matrix. Sample second moment matrices arise as 
central objects of interest in many areas, such as multivariate analysis, 
stochastic linear regression, time series analysis, and learning theory. 

Our motivation comes from learning theory and, in particular, from random 
sampling of trigonometric polynomials. Random sampling is a strategy of 
choice for learning an unknown function in a given class of functions. This 
idea is predominant in the version of learning theory and sampling theory by 
Cucker, Smale, and Zhou~\cite{smale2,SZ04}. In~\cite{bagr04} the
randomization of the samples was used for the justification of  
numerical algorithms. Random sampling and random measurements are central in 
the emerging field of sparse reconstruction, also referred to as compressed 
sensing \cite{carota06,cataXX,do04,kura06,ra05-7,ra06,ru06-1}. %

We first revisit the random sampling of trigonometric\ polynomials 
with a given degree or support, which was  studied in~\cite{bagr04}. We 
review and supplement the probability inequalities for the condition 
number of the associated Fourier sample second moment matrix in~\cite{bagr04} 
(Section \ref{Sec_2}). 

In  the main part (Section \ref{Sec_3})  we replace Fourier matrices by 
general random matrices with independent rows and derive probability 
estimates for sample second moment matrix obtained from general random 
matrices $U$. Our main result is an exponential probability inequality for 
the condition number of the sample second moment matrix for a vast class of 
random matrices.  The assumptions are extremely general:

(i) We treat random matrices with unbounded entries for which certain
moment generating functions exist  (Section~3.1)

(ii) We assume that the rows of the random matrix  are independent,
but we do not assume that the rows are identically distributed. This
somewhat technical  extension is treated in Section~\ref{noniid}. 

(iii) A further feature of our results is that all constants are given explicitly 
as a function of the parameters that describe the distribution of the random 
matrices. The explicit form of the constants is important to determine the 
sample size for which the condition number of the sample second moment 
matrix is small with high (``overwhelming'') probability. 

The proofs are based on  versions of Bernstein's inequality for sums
of (unbounded) random variables and suitable estimates for the
operator norm of a matrix. 

Under the standard assumptions for  random matrices, namely  bounded entries
and i.i.d.\ rows, several methods are available for estimates of the
sample second moment matrix. We mention in particular~\cite{BBZ,SZ07}
where exponential inequalities for operator-valued random variables
are used prominently. In \cite{Ru99} the non-commutative Khinchine
inequality is used to obtain an estimate for the expectation of the
operator norm of $n^{-1}U^*U - \mathrm{Id}$. 
A   beautiful,   and deep inequality for the sample second moment matrix
of random matrices was recently obtained by  Mendelson and
Pajor~\cite{mepaXX}.
Clearly, under significantly more restrictive assumptions on the
random matrix more precise exponential inequalities can be obtained. 
A detailed comparison between~\cite{mepaXX} and our result will be 
given in Section 3.1.1.

In learning theory one is often interested in the efficiency of the sampling 
procedure, i.e., in Cucker and Smale's words~\cite{smale2}, 
\textquotedblleft how many random samples do we need to assert, with 
confidence $1-\delta $, that the condition number does not exceed a given 
threshold.\textquotedblright\ For random sampling of trigonometric 
polynomials in several variables inspection shows that the probability 
inequalities in~\cite{bagr04}, as well as the ones in Section \ref{Sec_3} of 
the present paper, lead to lower bounds for the required sample size that 
are typically of the order $D^{2}\ln D$. We show in Section \ref{Sec_Comp} 
how the result in \cite{mepaXX} can be used to improve this order to $D\ln D$%
. In Section \ref{Sec_4} this result is further improved by using the method 
developed in~\cite{ra05-7} (after inspiration from~\cite{carota06}). To put 
it more casually, these results show that we need $const\cdot D\ln D$ random 
samples to learn a trigonometric\ polynomial taken from a $D$-dimensional 
space. This seems to be the optimal order that can be expected in a 
probabilistic setting. 

\textbf{Notation.} By $\Vert \cdot \Vert _{2}$ we denote the usual Euclidean 
norm on $\mathbb{C}^{D}$. For a (hermitian) matrix $A$ we denote by $\lambda 
_{\max }(A)$ and $\lambda _{\min }(A)$ the maximal and minimal eigenvalues 
of $A$. The condition number of $A$ is then given by $\lambda _{\max 
}(A)/\lambda _{\min }(A)$. For a matrix $A$ its transpose is denote by $%
A^{\prime }$ and its conjugate-transpose by $A^{\ast }$. The operator norm 
of a matrix is $\Vert A\Vert =\lambda _{\max }(A^{\ast }A)^{1/2}$. By $%
\mathbb{P}$ we denote the probability measure on the probability space 
supporting all the random variables used subsequently, and $\mathbb{E}$ 
denotes the corresponding expectation operator. 

\section{Random Sampling of Trigonometric Polynomials\label{Sec_2}} 

Let $\Gamma $ be a (non-empty) finite subset of $\mathbb{Z}^{d}$. By $%
\mathcal{P}_{\Gamma }$ we denote the space of all trigonometric polynomials 
in dimension $d$ with coefficients supported on $\Gamma $. Such a polynomial 
has the form  
\begin{equation*} 
f(x)\,=\,\sum_{k\in \Gamma }a_{k}e^{2\pi ik\cdot x},\quad x\in \lbrack 
0,1]^{d} 
\end{equation*}%
with coefficients $a_{k}\in \mathbb{C}$. If $\Gamma =\{-m,-m+1,\ldots 
,m-1,m\}^{d}$, then $\mathcal{P}_{\Gamma }$ is the space of all 
trigonometric polynomials of degree at most $m$. We let $D=|\Gamma |$ be the 
dimension of $\mathcal{P}_{\Gamma }$. 

Let $x_{1},\ldots ,x_{{n}}\in \lbrack 0,1]^{d}$. We are interested in the 
reconstruction of a trigonometric\ polynomial $f$ from its sample values $%
f(x_{1}),\ldots ,f(x_{{n}})$. Let $y=(f(x_{1}),\ldots ,f(x_{n}))^{\prime }$ 
be the vector of sampled values of $f$ and let $U$ be the ${n}\times {D}$ 
matrix with entries%
\begin{equation} 
u_{tk}=e^{2\pi ik\cdot x_{t}},\quad k\in \Gamma ,\text{ }t=1,\ldots ,{n}. 
\label{defU} 
\end{equation}%
The reconstruction of $f$ amounts to solving the linear system%
\begin{equation*} 
Ua=y 
\end{equation*}%
for the coefficient vector $a=(a_{k})_{k\in \Gamma }$. Alternatively, one 
may try to solve the normal equation%
\begin{equation*} 
U^{\ast }Ua=U^{\ast }y. 
\end{equation*}%
We note that the invertibility of $U^{\ast }U$ is equivalent to the sampling 
inequalities  
\begin{equation} 
A\Vert f\Vert _{2}^{2}\leq \sum_{t=1}^{n}|f(x_{t})|^{2}=a^{\ast }U^{\ast 
}Ua\leq B\Vert f\Vert _{2}^{2}\qquad \text{for all}\,f\in \mathcal{P}%
_{\Gamma }\,,  \label{eq:c3} 
\end{equation}%
for some positive real numbers $A$ and $B$, and that the condition number of  
$U^{\ast }U$ is bounded by $B/A$. 

In the spirit of learning theory, one assumes that the sampling points are 
taken at random. Then the matrix $U^{\ast }U$ is a random matrix depending 
on the sampling points $(x_{t})$. Several questions arise: 

\begin{enumerate} 
\item Determine the probability that $U^{\ast }U$ is invertible. 

\item Determine the probability that the condition number of $U^{\ast }U$ 
does not exceed a given threshold. 

\item Determine the number of random samples required to achieve such 
estimates. This is the effectivity problem for random sampling. 
\end{enumerate} 

For trigonometric\ polynomials Question 1 has been answered in~\cite[Thm.~1.1%
]{bagr04}: If the $x_{t}$ are i.i.d.~with a distribution that is absolutely 
continuous with respect to the Lebesgue-measure on $[0,1]^{d}$, then $%
U^{\ast }U$ is invertible almost surely provided ${n}\geq {D}$. 
Some answers to Questions 2 and 3 are also provided in~\cite%
{bagr04}. It is shown that $U^{\ast }U$ is well-conditioned whenever the 
number of samples is large enough~\cite[Thms.~5.1, 6.2]{bagr04}: 

\begin{Theorem} 
\label{thm:BassGr} Assume that $x_{1},\ldots ,x_{n}$ are i.i.d.~random 
variables uniformly distributed on $[0,1]^{d}$. Let $U$ be the associated 
random Fourier matrix defined in (\ref{defU}). Let $\varepsilon \in (0,1)$. 
There exist positive constants $A,B$ depending only on ${D}=|\Gamma |$ such 
that the event  
\begin{equation*} 
1-\varepsilon \leq \lambda _{\min }({n}^{-1}U^{\ast }U)\leq \lambda _{\max }(%
{n}^{-1}U^{\ast }U)\leq 1+\varepsilon 
\end{equation*}%
has probability at least%
\begin{equation} 
1-Ae^{-B{n}\varepsilon ^{2}/(1+\varepsilon )}.  \label{eq:c1} 
\end{equation}%
In particular, with probability not less than (\ref{eq:c1}), the condition 
number of $U^{\ast }U$ is bounded by $(1+\varepsilon )/(1-\varepsilon )$. 
\end{Theorem} 

A careful analysis of the constants $A$ and $B$ in~(\ref{eq:c1}) reveals 
that the number of samples required in (\ref{eq:c1}) to guarantee a 
probability $\geq 1-\delta $ is%
\begin{equation} 
{n}\,\geq \,C{D}^{2}\ln {D},  \label{n:min1} 
\end{equation}%
where $C$ depends on $\delta $ and $\varepsilon $. If $\Gamma =\{-m,\dots 
,m\}^{d}$ (trigonometric\ polynomials of degree $m$ in $d$ variables), then $%
D=(2m+1)^{d}$, and this bound on the number of samples is unfortunately too 
large to be useful at more realistic sample sizes like $n\sim D$ or $n\sim 
D\ln D$. 

For the case $\Gamma =\{-m,\dots ,m\}^{d}$  
a better estimate for the condition number can be extracted from~\cite%
{bagr04}. We work, however, with a slightly different matrix. Given the 
sampling points $x_{1},\ldots ,x_{n}$, we define the Voronoi regions  
\begin{equation*} 
V_{t}\,:=\,\{y\in \lbrack 0,1]^{d}:\Vert y-x_{t}\Vert _{2}\leq \Vert 
y-x_{s}\Vert _{2},s\neq t,1\leq s\leq {n}\},\quad t=1,\ldots ,{n} 
\end{equation*}%
and let $w_{t}=|V_{t}|$ be the Lebesgue measure of $V_{t}$. We consider the 
weighted matrix $T^{w}$  
\begin{equation*} 
T^{w}\,:=\,U^{\ast }WU 
\end{equation*}%
where $W$ is the diagonal matrix with the weights $w_{t}$, $t=1,\ldots ,{n}$%
, on the diagonal. Note that $a$ is also the solution of $T^{w}a=U^{\ast }Wy$%
. The following result is implicit in~\cite{bagr04}. 

\begin{Theorem} 
\label{thm:detprob} Let $\Gamma =\{-m,\ldots ,m\}^{d}$, i.e., we consider 
trigonometric polynomials of $d$ variables of degree $m$. Suppose that $%
x_{1},\ldots ,x_{n}$ are i.i.d~random variables which are uniformly 
distributed on $[0,1]^{d}$. Choose $\gamma \in (0,1)$. If  
\begin{equation} 
{n}\geq \left( \frac{2\pi d}{\gamma \ln 2}\right) ^{d}m^{d}\ln \left( \left(  
\frac{2\pi d}{\gamma \ln 2}\right) ^{d} \, \frac{m^{d}}{\delta } \right) \,, 
\label{eq:c2} 
\end{equation}%
then with probability at least $1-\delta $ the condition number of $T^{w}$ 
is bounded by $\left( 1-2^{\gamma -1}\right) ^{-2}$. 
\end{Theorem} 

\begin{Proof} 
By combining a deterministic estimate with a probabilistic covering result, 
the following estimate was derived in~\cite[Thm.~4.2]{bagr04}: Let $N\in  
\mathbb{N}$ be arbitrary; then with probability at least $%
1-N^{d}e^{-n/N^{d}} $ we have  
\begin{equation*} 
(2-e^{2\pi md/N})^{2}\leq \lambda _{\min }(T^{w})\leq \lambda _{\max 
}(T^{w})\leq 4\,. 
\end{equation*}%
For the condition number to be bounded by $4\left( 2-2^{\gamma }\right) 
^{-2} $ with probability at least $1-\delta $, we need that  
\begin{equation*} 
2\pi md/N\leq \gamma \ln 2\qquad \text{and}\quad N^{d}e^{-n/N^{d}}\leq 
\delta \,. 
\end{equation*}%
By solving for $n$, we find that $n$ must satisfy the inequality~%
\eqref{eq:c2}. 
\end{Proof} 

Since ${D}=|\Gamma |=|\{-m,-m+1,\ldots ,m-1,m\}^{d}|=(2m+1)^{d}$, Theorem~%
\ref{thm:detprob} becomes effective for  
\begin{equation} 
n\approx (\pi d)^{d}{D}\ln \left( (\pi d)^{d}D\right) \,. 
\label{est:detprob} 
\end{equation}%
Thus Theorem~\ref{thm:detprob} is a genuine improvement over (\ref{n:min1}) 
for fixed value of $d$. The dependence $n\sim D\ln D$ on the dimension of 
the function space seems to be of the correct order. However, the constant $%
\left( \pi d/\gamma \right) ^{d}$ depends strongly on the number of 
variables $d$, and so Theorem~\ref{thm:detprob} does not escape the curse of 
dimensionality.

In Theorem \ref{thm_uniform} we will prove a much better result for the 
condition number of $U^{\ast }U$ where the constants do neither depend on $d$ 
nor on the special form of the spectrum $\Gamma $. See also Corollary \ref%
{imprv}. 

\section{Exponential Inequalities For Sample Second Moment Matrices \label%
{Sec_3}} 

In this section we abstract from the concrete form of $U$ as given in (\ref%
{defU}) and consider arbitrary complex random matrices with independent rows 
satisfying some regularity conditions. Apart from being of interest in its 
own, this more general setting allows one to study random sampling not only 
for trigonometric polynomials but also for more general types of 
finite-dimensional function spaces, such as random sampling of algebraic 
polynomials on domains, or of spaces of spherical harmonics on the sphere 
(see~\cite[Sect.~6]{bagr04} for a list of examples). 

For the sake of exposition and clarity we first treat the case
of i.i.d.~rows and then only the general case.

\subsection{The I.I.D. Case\label{iid}} 

We assume first that the random matrix $U\in \mathbb{C}^{n\times D}$ has 
independent identically distributed rows and delay the discussion of the 
case of independent, but not identically distributed rows to Section \ref%
{noniid}. Furthermore, we assume that the rows $u_{t\cdot }=(u_{t1},\ldots 
,u_{tD})$ of $U$ satisfy the following condition: The moment generating 
functions of the random variables $\Re (\overline{u_{1k}}u_{1j})$ and $\Im (%
\overline{u_{1k}}u_{1j})$ exist for all $1\leq k,j\leq D$; i.e., there 
exists $x_{0}>0$ such that for all $1\leq k,j\leq D$  
\begin{equation} 
\mathbb{E}\left[ \exp (x\Re (\overline{u_{1k}}u_{1j}))\right] <\infty 
,\qquad \mathbb{E}\left[ \exp (x\Im (\overline{u_{1k}}u_{1j}))\right] <\infty 
\label{moment_gen_func_1} 
\end{equation}%
hold for all $x<x_{0}$. Note that a sufficient condition for (\ref%
{moment_gen_func_1}) is that the moment generating function of $\left\vert 
u_{1k}\right\vert ^{2}+\left\vert u_{1j}\right\vert ^{2}$ exists for all $%
k,j $. Further, we let%
\begin{equation*} 
Q:=\mathbb{E(}u_{1\cdot }^{\ast }u_{1\cdot })\in \mathbb{C}^{{D}\times {D}} 
\end{equation*}%
with entries $q_{kj}$. We note that by the strong law of large numbers $%
n^{-1}U^{\ast }U$ converges to $Q=\mathbb{E}[{n}^{-1}U^{\ast }U]$ almost 
surely. 

Assumption (\ref{moment_gen_func_1}) is easily seen to be equivalent to the 
existence of finite constants $M\geq 0$ and $v_{kj}\geq 0$ such that for all  
$\ell \geq 2$  
\begin{align} 
\mathbb{E}\left[ |\Re (\overline{u_{1k}}u_{1j}-q_{kj})|^{\ell }\right] \,& 
\leq \,2^{-1}\ell !\,M^{\ell -2}v_{kj},  \label{moment_estimateRe_1} \\ 
\mathbb{E}\left[ |\Im (\overline{u_{1k}}u_{1j}-q_{kj})|^{\ell }\right] \,& 
\leq \,2^{-1}\ell !\,M^{\ell -2}v_{kj}  \label{moment_estimateIm_1} 
\end{align}%
hold for all $1\leq k,j\leq D$. For a generalization leading to a slightly 
better, but more complex bound see Section \ref{noniid}. 

\begin{remark} 
\label{rem:bounded} If the random variables $u_{1k}$ are bounded, i.e.,  
\begin{equation*} 
|\Re (\overline{u_{1k}}u_{1j}-q_{kj})|\,\leq \,C\quad \text{and}\quad |\Im (%
\overline{u_{1k}}u_{1j}-q_{kj})|\,\leq \,C 
\end{equation*}%
holds with probability $1$ for all $1\leq k,j\leq D$, then (\ref%
{moment_estimateRe_1}) and (\ref{moment_estimateIm_1}) hold with $M=C/3$ and  
\begin{equation} 
v_{kj}=\max \{\mathbb{E}\left[ (\Re (\overline{u_{1k}}u_{1j}-q_{kj}))^{2}%
\right] ,\mathbb{E}\left[ (\Im (\overline{u_{1k}}u_{1j}-q_{kj}))^{2}\right] 
\}.  \label{vij_bounded_1} 
\end{equation}%
This claim is obvious for $\ell =2$. For $\ell \geq 3$ it follows from a 
general inequality for arbitrary real-valued bounded random variables $X$:  
\begin{equation*} 
\mathbb{E}\left[ |X|^{\ell }\right] \,=\,\mathbb{E}\left[ X^{2}|X|^{\ell -2}%
\right] \,\leq \,C^{\ell -2}\mathbb{E}[X^{2}]\,=\,\ell !\,\frac{C^{\ell -2}}{%
\ell !}\sigma ^{2}\,\leq \,\ell !\frac{C^{\ell -2}}{2\cdot 3^{\ell -2}}%
\sigma ^{2}, 
\end{equation*}%
where $\sigma ^{2}=\mathbb{E}[X^{2}]$ and $|X|\leq C$ holds with probability  
$1$. In particular, this shows that the random Fourier matrix given in (\ref%
{defU}) satisfies (\ref{moment_estimateRe_1}) and (\ref{moment_estimateIm_1}%
).$\hfill \square $ 
\end{remark} 

The proof of the main result in this section will make use of the following 
Bernstein-type inequality for unbounded random variables given in Bennett  
\cite[eq.~(7)]{be62}, see also \cite[Lemma~2.2.11]{vawe96}: 

\textit{Let }$X_{1},\ldots ,X_{{n}}$\textit{\ be independent real-valued 
random variables with zero mean such that }$E|X_{t}|^{\ell }\leq \ell 
!M^{\ell -2}v_{t}/2$\textit{\ holds for every }$\ell \geq 2$\textit{\ and }$%
t=1,\ldots ,n$\textit{\ for some finite constants }$M\geq 0$\textit{\ and }$%
v_{t}\geq 0$\textit{. Then for every }$x>0$\textit{\ }%
\begin{equation} 
\mathbb{P}\left( \left\vert \sum_{t=1}^{{n}}X_{t}\right\vert \geq x\right) 
\leq 2e^{-\frac{x^{2}}{2}\left( \sum_{t=1}^{{n}}v_{t}+Mx\right) ^{-1}}, 
\label{bennett} 
\end{equation}%
\textit{with the convention that the right-hand side in (\ref{bennett}) is 
zero if }$M=0$\textit{\ and }$\sum_{t=1}^{{n}}v_{t}=0$\textit{.} 

Note that Bennett \cite{be62} assumes $\sum_{t=1}^{{n}}v_{t}>0$ but the 
inequality (\ref{bennett}) trivially also holds for $\sum_{t=1}^{{n}}v_{t}=0$ 
in which case the probability on the left-hand side is zero. Inequality (\ref%
{bennett}), and hence the subsequent results, can be somewhat improved, see 
Bennett \cite[eq.~(7a)]{be62}. Since this does not result in any significant 
gain, we do not give the details. 

Set  
\begin{equation*} 
v:=\max_{1\leq k,j\leq D}v_{kj}. 
\end{equation*}%
Note that neither $v$ nor $M$ depend on $n$ because the rows are identically 
distributed. However, they depend on the distribution of the random vector $%
u_{1\cdot }$ and hence may depend on $D$. Our main result now reads as 
follows. 

\begin{Theorem} 
\label{thm:main} Assume that the rows $u_{1\cdot },\ldots ,u_{{n\cdot }}$ of  
$U$ are 
i.i.d.~random vectors in $\mathbb{C}^{{D}}$ whose entries satisfy the moment 
bounds (\ref{moment_estimateRe_1}) and (\ref{moment_estimateIm_1}). Then, 
for every $\varepsilon >0$, the operator norm satisfies  
\begin{equation*} 
\left\Vert {n}^{-1}U^{\ast }U-Q\right\Vert <\varepsilon 
\end{equation*}%
with probability at least%
\begin{equation} 
1-4{D}^{2}\exp \left( -\frac{n\varepsilon ^{2}}{{D}^{2}\left( 4v+2\sqrt{2}%
D^{-1}M\varepsilon \right) }\right) .  \label{prob_bound_1} 
\end{equation}%
In particular, with probability not less than (\ref{prob_bound_1}) the 
extremal eigenvalues of $n^{-1}U^{\ast }U$ satisfy%
\begin{equation} 
\lambda _{\min }(Q)-\varepsilon <\lambda _{\min }({n}^{-1}U^{\ast }U)\leq 
\lambda _{\max }({n}^{-1}U^{\ast }U)<\lambda _{\max }(Q)+\varepsilon \, . 
\label{eig_estim_1} 
\end{equation}%
Consequently, if $Q$ is non-singular and $\varepsilon \in (0,\lambda _{\min 
}(Q))$, then the condition number of $U^{\ast }U$ is bounded by $\frac{%
\lambda _{\max }(Q)+\varepsilon }{\lambda _{\min }(Q)-\varepsilon }$ with 
probability not less than (\ref{prob_bound_1}). 
\end{Theorem} 

In connection with (\ref{eig_estim_1}) we note that $\lambda _{\min }({n}%
^{-1}U^{\ast }U)\geq 0$ holds trivially, since the matrix ${n}^{-1}U^{\ast 
}U $ is nonnegative definite. 

\bigskip 

\begin{Proof} 
We first note that inequality (\ref{eig_estim_1}) for the extremal 
eigenvalues of $n^{-1}U^{\ast }U$ follows from the inequality $\left\Vert {n}%
^{-1}U^{\ast }U-Q\right\Vert <\varepsilon $ for the operator norm. Hence, it 
suffices to concentrate on the operator norm, which we majorize with Schur's 
test by using that $\Vert A\Vert \leq \max_{k}\sum_{j}|a_{kj}|$ for 
self-adjoint $A$. In this way we obtain that  
\begin{align} 
& \mathbb{P}(\left\Vert {n}^{-1}U^{\ast }U-Q\right\Vert \geq \varepsilon 
)\leq \mathbb{P}\left( \max_{k=1,...,{D}}\sum_{j=1}^{{D}}\left\vert {n}%
^{-1}\sum_{t=1}^{{n}}(\overline{u_{tk}}u_{tj}-q_{kj})\right\vert \geq 
\varepsilon \right)  \notag \\ 
& \leq \sum_{k=1}^{{D}}\mathbb{P}\left( \sum_{j=1}^{{D}}\left\vert {n}%
^{-1}\sum_{t=1}^{{n}}(\overline{u_{tk}}u_{tj}-q_{kj})\right\vert \geq 
\varepsilon \right) \leq \sum_{k,j=1}^{{D}}\mathbb{P}\left( \left\vert {n}%
^{-1}\sum_{t=1}^{{n}}(\overline{u_{tk}}u_{tj}-q_{kj})\right\vert \geq 
\varepsilon /{D}\right)  \notag \\ 
& =\sum_{k,j=1}^{{D}}\mathbb{P}\left( \left\vert {n}^{-1}\sum_{t=1}^{{n}}(%
\overline{u_{tk}}u_{tj}-q_{kj})\right\vert ^{2}\geq (\varepsilon /{D}%
)^{2}\right)  \label{Schur} \\ 
& =\sum_{k,j=1}^{{D}}\mathbb{P}\left( \left\vert {n}^{-1}\sum_{t=1}^{{n}}\Re 
(\overline{u_{tk}}u_{tj}-q_{kj})\right\vert ^{2}+\left\vert {n}%
^{-1}\sum_{t=1}^{{n}}\Im (\overline{u_{tk}}u_{tj}-q_{kj})\right\vert 
^{2}\geq (\varepsilon /{D})^{2}\right)  \notag \\ 
& \leq \sum_{k,j=1}^{{D}}\mathbb{P}\left( \left\vert \sum_{t=1}^{{n}}\Re (%
\overline{u_{tk}}u_{tj}-q_{kj})\right\vert \geq \frac{{n}\varepsilon }{\sqrt{%
2}{D}}\right) +\sum_{k,j=1}^{{D}}\mathbb{P}\left( \left\vert \sum_{t=1}^{{n}%
}\Im (\overline{u_{tk}}u_{tj}-q_{kj})\right\vert \geq \frac{{n}\varepsilon }{%
\sqrt{2}{D}}\right) .  \notag 
\end{align}%
For each index $k,j$ the inequality (\ref{bennett}) gives  
\begin{equation} 
\mathbb{P}\left( \left\vert \sum_{t=1}^{{n}}\Re (\overline{u_{tk}}%
u_{tj}-q_{kj})\right\vert \geq \frac{{n}\varepsilon }{\sqrt{2}{D}}\right) 
\leq 2\exp \left( -\frac{n\varepsilon ^{2}}{{D}^{2}\left( 4v_{kj}+2\sqrt{2}%
D^{-1}M\varepsilon \right) }\right)  \label{ben} 
\end{equation}%
and similarly for the imaginary part. Hence, we finally obtain  
\begin{eqnarray} 
&&\mathbb{P}\left( \Vert {n}^{-1}U^{\ast }U-Q\Vert \geq \varepsilon \right) 
\,\leq 4\sum_{k,j=1}^{{D}}\exp \left( -\frac{n\varepsilon ^{2}}{{D}%
^{2}\left( 4v_{kj}+2\sqrt{2}D^{-1}M\varepsilon \right) }\right)  \notag \\ 
&\leq &\,4{D}^{2}\exp \left( -\frac{n\varepsilon ^{2}}{{D}^{2}\left( 4v+2%
\sqrt{2}D^{-1}M\varepsilon \right) }\right)  \label{union_bound} 
\end{eqnarray}%
with $v$ as defined above. Thus (\ref{prob_bound_1}) follows. 
\end{Proof} 

\begin{remark} 
\label{nonsing}If $u_{1.}$ possesses an absolutely continuous distribution 
then $Q$ is automatically non-singular. More generally, this holds as long 
as the distribution of $u_{1.}$ is not concentrated on a $(D-1)$-dimensional 
linear subspace of $\mathbb{C}^{D}$. To see this, consider the quadratic 
forms $z^{\ast }Qz$ for $z\in \mathbb{C}^{D}$ and note that%
\begin{equation*} 
z^{\ast }Qz=\mathbb{E}[z^{\ast }u_{1\cdot }^{\ast }u_{1\cdot }z]=\mathbb{E}%
[\left\vert u_{1\cdot }z\right\vert ^{2}]\geq 0. 
\end{equation*}%
Hence, if $z^{\ast }Qz=0$, then $\left\vert u_{1\cdot }z\right\vert ^{2}=0$ 
with probability $1$; thus the distribution of $u_{1.}$ would have to reside 
in the orthogonal complement of the one-dimensional subspace spanned by $%
z^{\ast }$.$\hfill \square $ 
\end{remark} 

\begin{remark} 
\label{real}For real-valued random matrices $U$ we can improve the 
probability bound (\ref{prob_bound_1}) to  
\begin{equation*} 
1-2{D}^{2}\exp \left( -\frac{{n}\varepsilon ^{2}}{2{D}^{2}\left( v+\frac{%
M\varepsilon }{{D}}\right) }\right) . 
\end{equation*}%
A similar improvement for real-valued $U$ applies to the subsequent 
corollary and remark as well as to the results in Section \ref{noniid}.$%
\hfill \square $ 
\end{remark} 

\begin{corollary} 
\label{cor:bounded} Assume that the rows $u_{1\cdot },\ldots ,u_{{n\cdot }}$ 
of $U$ are 
i.i.d.~random vectors in $\mathbb{C}^{{D}}$ that are bounded, i.e.,  
\begin{equation*} 
|\Re (\overline{u_{1k}}u_{1j}-q_{kj})|\,\leq \,C\quad and\quad |\Im (%
\overline{u_{1k}}u_{1j}-q_{kj})|\,\leq \,C 
\end{equation*}%
holds with probability $1$ for or all $1\leq k,j\leq D$. Let  
\begin{equation*} 
b\,:=\,\max_{k,j=1,\ldots ,{D}}\left\{ \mathbb{E}\left[ (\Re (\overline{%
u_{1k}}u_{1j}-q_{kj}))^{2}\right] ,\mathbb{E}\left[ (\Im (\overline{u_{1k}}%
u_{1j}-q_{kj}))^{2}\right] \right\} . 
\end{equation*}%
Then the conclusions of Theorem \ref{thm:main} hold and (\ref{prob_bound_1}) 
becomes%
\begin{equation} 
1-4{D}^{2}\exp \left( -\frac{{n}\varepsilon ^{2}}{{D}^{2}\left( 4b+2\sqrt{2}%
D^{-1}C\varepsilon /3\right) }\right) .  \label{6} 
\end{equation} 
\end{corollary} 

\begin{Proof} 
By Remark \ref{rem:bounded} conditions (\ref{moment_estimateRe_1}) and (\ref%
{moment_estimateIm_1}) hold with $M=C/3$ and $v_{kj}$ as in (\ref%
{vij_bounded_1}). Then the statement follows from Theorem \ref{thm:main}. 
\end{Proof} 

\begin{remark} 
\label{rem:bounded_1}Corollary~\ref{cor:bounded} can also be derived by 
using the classical Bernstein inequality instead of inequality (\ref{bennett}%
) in the proof of Theorem \ref{thm:main}. Furthermore, the bound in (\ref{6}%
) can be somewhat improved by using an improved form of Bernstein's 
inequality \cite[eq.~(8)]{be62} (see also \cite[Corollary A.2]{celu06}) for 
bounded random variables instead of (\ref{bennett}) in that step: If we use 
that inequality in the estimate (\ref{ben}), we arrive at the following 
improved bound (provided $C>0$, $b>0$):%
\begin{equation*} 
1-4{D}^{2}\exp \left( -C^{-2}{n}b\left( \left( 1+\frac{C\varepsilon }{\sqrt{2%
}{D}b}\right) \ln \left( 1+\frac{C\varepsilon }{\sqrt{2}{D}b}\right) -\frac{%
C\varepsilon }{\sqrt{2}{D}b}\right) \right) . 
\end{equation*}%
$\hfill \square $ 
\end{remark} 

Let us now apply our findings to random sampling of trigonometric 
polynomials. 

\begin{corollary} 
Let $x_{1},\ldots ,x_{{n}}$ be independent random variables uniformly 
distributed on $[0,1]^{d}$. Let $U$ be the associated ${n}\times {D}$ random 
Fourier matrix (\ref{defU}). Let $\varepsilon >0$. Then with probability at 
least  
\begin{equation} 
1-4{D}({D}-1)\exp \left( -\frac{{n}\varepsilon ^{2}}{2\left( ({D-1)}^{2}+%
\sqrt{2}(D-1)\varepsilon /3\right) }\right)  \label{pprob} 
\end{equation}%
we have%
\begin{equation} 
\left\Vert {n}^{-1}U^{\ast }U-Q\right\Vert <\varepsilon , 
\label{operatornorm} 
\end{equation}%
and hence  
\begin{equation} 
1-\varepsilon <\lambda _{\min }({n}^{-1}U^{\ast }U)\leq \lambda _{\max }({n}%
^{-1}U^{\ast }U)<1+\varepsilon \,.  \label{Fourier:eigs} 
\end{equation}%
Consequently, for $0<\varepsilon <1$, the condition number of $U^{\ast }U$ 
is bounded by $(1+\varepsilon )/(1-\varepsilon )$ with probability not less 
than (\ref{pprob}). 
\end{corollary} 

\begin{Proof} 
In this case $({n}^{-1}U^{\ast }U)_{kj}={n}^{-1}\sum_{t=1}^{{n}}e^{2\pi 
i(j-k)\cdot x_{t}}$ and consequently $Q=I$, so $\lambda _{\min }(Q)=\lambda 
_{\max }(Q)=1$. [By abuse of notation, $k$ denotes both an element of $%
\Gamma $ and a column index.] Furthermore, $\sum_{t=1}^{{n}}(\overline{u_{tk}%
}u_{tj}-q_{kj})=0$ for $k=j$. Hence, the double sum in the second line of (%
\ref{Schur}) only extends over $j\neq k$ and consequently $\varepsilon /D$ 
can be replaced by $\varepsilon /(D-1)$ in the subsequent steps in (\ref%
{Schur}). Furthermore, when deducing (\ref{union_bound}) from the union 
bound in (\ref{Schur}) we only have to take into account ${D}({D}-1)$ 
instead of ${D}^{2}$ summands; cf. also Remark \ref{zero} below. Moreover, $%
|\Re (\overline{u_{tk}}u_{tj}-q_{kj})|\leq 1$ for all $k,j$. For $k\neq j$ 
we have  
\begin{equation*} 
\mathbb{E}\left[ \Re (\overline{u_{tk}}u_{tj}-q_{kj})^{2}\right] 
\,=\,\int_{[0,1]^{d}}(\Re (\exp (2\pi i(j-k)\cdot x))^{2}dx\,=\,\frac{1}{2}, 
\end{equation*}%
hence $v_{kj}=1/2$. The same holds for the imaginary part. In view of Remark %
\ref{rem:bounded} the result follows. 
\end{Proof} 

From the previous result it is easy to determine the minimal number of 
sampling points sufficient to provide a small condition number with high 
probability. 

\begin{corollary} 
\label{Fourier}Let $x_{1},\ldots ,x_{{n}}$ be independent random variables 
uniformly distributed on $[0,1]^{d}$. Let $U$ be the associated ${n}\times {D%
}$\ random Fourier matrix (\ref{defU}). Let $0<\varepsilon <1,0<\delta <1$ 
and suppose  
\begin{equation} 
{n}\geq \frac{2}{\varepsilon ^{2}}\left( {(D-1)}^{2}+\frac{\sqrt{2}({D-1)}%
\varepsilon }{3}\right) \ln \left( \frac{4{D}({D}-1)}{\delta }\right) . 
\label{minimum:n} 
\end{equation}%
Then (\ref{operatornorm}) and (\ref{Fourier:eigs}) hold with probability at 
least $1-\delta $. 
\end{corollary} 

We note that (\ref{minimum:n}) is implied by the more compact inequality  
\begin{equation} 
{n}\geq \frac{C{D}^{2}\ln ({D}/\delta )}{\varepsilon ^{2}} 
\label{minimum:nb} 
\end{equation}%
for an appropriate constant $C$. We will improve on this result in Corollary %
\ref{imprv} and in Section \ref{Sec_4}, see in particular (\ref{minimum:n2}). 

\subsubsection{Comparison With Other Results\label{Sec_Comp}} 

Recently Mendelson and Pajor \cite{mepaXX} provide a related exponential 
inequality for 
random matrices with i.i.d. real-valued rows. They assume the following 
properties: 

(a) There exists $\rho >0$ such that for every $\theta \in \mathbb{R}%
^{D},\Vert \theta \Vert _{2}=1$, $(\mathbb{E}|\langle u_{1\cdot },\theta 
\rangle |^{4})^{1/4}\leq \rho <\infty $. 

(b) Set $Z=\Vert u_{1\cdot }\Vert _{2}$, then $\Vert Z\Vert _{\psi _{\alpha 
}}<\infty $ for some $\alpha \geq 1$. 

Here the Orlicz norm $\Vert \cdot \Vert _{\psi _{\alpha }}$ of a real-valued 
random variable $Y$ with respect to $\psi _{\alpha }(x)=\exp (x^{\alpha })-1$ 
is defined as $\Vert Y\Vert _{\psi _{\alpha }}\,=\,\inf \left\{ C>0:\mathbb{E%
}\psi _{\alpha }\left( |Y|/C\right) \leq 1\right\} $. Note that if $\alpha 
\geq 2$, 
condition (b) is stronger than our assumption (\ref{moment_gen_func_1}), 
hence condition (b) implies (\ref{moment_estimateRe_1})-(\ref%
{moment_estimateIm_1}). 

Under conditions (a) and (b), Mendelson and Pajor \cite[Theorem 2.1]{mepaXX}%
\ show that there exists an absolute constant $c>0$ such that for every $%
\varepsilon >0$ the operator norm satisfies%
\begin{equation} 
\mathbb{P(}\left\Vert {n}^{-1}U^{\ast }U-Q\right\Vert <\varepsilon )\geq 
1-2\exp \left( -\left( \frac{c\varepsilon }{\max \{B_{{n}},A_{{n}}^{2}\}}%
\right) ^{\alpha /(\alpha +2)}\right) \, ,  \label{MP} 
\end{equation}%
where%
\begin{equation*} 
A_{{n}}\,=\,\Vert Z\Vert _{\psi _{\alpha }}\frac{\sqrt{\ln (\min (D,n))}(\ln  
{n})^{1/\alpha }}{\sqrt{{n}}},\quad \quad B_{{n}}\,=\,\frac{\rho ^{2}}{\sqrt{%
{n}}}+\lambda _{\max }^{1/2}(Q)A_{{n}}. 
\end{equation*} 

We have added a factor $2$ on the right-hand side of (\ref{MP}) to correct a 
missing constant in~\cite[Theorem 2.1]{mepaXX}.  
Since the constant $c$ is not specified, the value of \eqref{MP} is mainly 
for asymptotics as $n\rightarrow \infty $, whereas the results of Section 3 
yield estimates with explicit constants for given $n$. Moreover, the 
probability estimate \eqref{MP} is only subexponential. For fixed dimension $%
D$, the right-hand side of (\ref{MP}) is of the order  
\begin{equation*} 
1-2\exp \left( -c_{1}n^{\alpha /(2\alpha +4)}(\ln n)^{-1/(\alpha +2)}\right) 
, 
\end{equation*}%
which is only subexponential, whereas the bound in (\ref{prob_bound_1}) is 
exponential of the form%
\begin{equation*} 
1-c_{2}\exp \left( -c_{3}n\right) . 
\end{equation*}%
Here $c_{1}$, $c_{2}$, and $c_{3}$ are constants that depend on $D$.

On the other hand, estimate~\eqref{MP} behaves much better with respect
to the dimension $D$. Indeed, \eqref{MP} can be used to improve
upon (\ref{minimum:n}) in the special case where the  
set $\Gamma $ is symmetric in the sense that $k\in \Gamma $ implies $-k\in 
\Gamma $. 


\begin{corollary} 
\label{imprv}Let $x_{1},\ldots ,x_{{n}}$ be independent random variables 
uniformly distributed on $[0,1]^{d}$. Let $U$ be the associated ${n}\times {D%
}$\ random Fourier matrix (\ref{defU}) and assume that $\Gamma $ is 
symmetric. 
Let $0<\varepsilon <1,0<\delta <1$ and suppose  
\begin{equation} 
{n}\geq \max \left\{ D,c^{-1}\varepsilon ^{-1}\ln \left( \frac{2}{\delta }%
\right) D\ln D,\left[ c^{-1}\varepsilon ^{-1}\ln \left( \frac{2}{\delta }%
\right) \right] ^{2}(\sqrt{D}+\sqrt{D\ln D})^{2}\right\} 
\label{better_bound} 
\end{equation}%
where $c$ is the absolute constant in (\ref{MP}). Then (\ref{operatornorm}) 
and (\ref{Fourier:eigs}) hold with probability at least $1-\delta $. 
\end{corollary} 

\begin{Proof} 
Although it is possible to adapt \cite{mepaXX} to complex-valued random 
matrices, we will use the result as stated. 

Since $\Gamma $ is symmetric by assumption, we may write it as $\Delta \cup 
(-\Delta )$ with $\Delta \cap (-\Delta )\subseteq \{0\}$. Define the real $%
n\times D$ matrix $W$ by $w_{t(-k)}=\sqrt{2}\cos (2\pi k\cdot x_{t})$ and $%
w_{tk}=\sqrt{2}\sin (2\pi k\cdot x_{t})$ for $k\in \Delta \backslash \{0\}$ 
and set $w_{t0}=1$ if $0\in \Delta $. Then clearly $U=WS$ where $S$ is a 
unitary $D\times D$ matrix and consequently $\left\Vert {n}^{-1}U^{\ast 
}U-I\right\Vert =\left\Vert {n}^{-1}W^{\ast }W-I\right\Vert $. 
To apply \eqref{MP} we note that $\Vert w_{1\cdot }\Vert _{2}=D^{1/2}$ and 
hence $\Vert Z\Vert _{\psi _{\alpha }}=\Vert \Vert w_{1\cdot }\Vert 
_{2}\Vert _{\psi _{\alpha }}=D^{1/2}(\ln 2)^{-1/\alpha }$ for every $\alpha 
\geq 1$. Furthermore,  
\begin{eqnarray*} 
\sup_{\Vert \theta \Vert _{2}=1,\theta \in \mathbb{R}^{D}}(\mathbb{E}%
|\langle w_{1\cdot },\theta \rangle |^{4})^{1/4} &=&\sup_{\Vert \theta \Vert 
_{2}=1,\theta \in \mathbb{R}^{D}}(\mathbb{E}|\langle u_{1\cdot }S^{\ast 
},\theta \rangle |^{4})^{1/4} \\ 
&=&\sup_{\Vert \theta \Vert _{2}=1,\theta \in \mathbb{R}^{D}}(\mathbb{E}%
|\langle u_{1\cdot },S\theta \rangle |^{4})^{1/4} \\ 
&\leq &\sup_{\Vert \theta \Vert _{2}=1,\theta \in \mathbb{C}^{D}}(\mathbb{E}%
|\langle u_{1\cdot },\theta \rangle |^{4})^{1/4}. 
\end{eqnarray*}%
Now, since $\left\vert \sum_{k\in \Gamma }\overline{\theta }_{k}\exp (2\pi 
ik\cdot x_{1})\right\vert ^{2}\leq \left( \sum_{k\in \Gamma }\left\vert 
\theta _{k}\right\vert \right) ^{2}$ we obtain  
\begin{eqnarray*} 
\mathbb{E}|\langle u_{1\cdot },\theta \rangle |^{4} &=&\mathbb{E}\left[ 
\left\vert \sum_{k\in \Gamma }\overline{\theta }_{k}\exp (2\pi ik\cdot 
x_{1})\right\vert ^{4}\right]  \\ 
&\leq &\left( \sum_{k\in \Gamma }\left\vert \theta _{k}\right\vert \right) 
^{2}\mathbb{E}\left[ \left\vert \sum_{k\in \Gamma }\overline{\theta }%
_{k}\exp (2\pi ik\cdot x_{1})\right\vert ^{2}\right]  \\ 
&\leq &D\Vert \theta \Vert _{2}^{2}=D. 
\end{eqnarray*}%
This shows that the rows $w_{t\cdot }$ satisfy condition (a) in \cite[%
Theorem 2.1]{mepaXX} with $\rho =D^{1/4}$. As a consequence, (\ref{MP}) 
applies to $W$, and hence to the Fourier matrix $U$, for every $\alpha \geq 1 
$. Since the left-hand side of (\ref{MP}) does not depend on $\alpha $, we 
may let $\alpha \rightarrow \infty $ and obtain for $n\geq D$ the bound%
\begin{equation*} 
\mathbb{P(}\left\Vert {n}^{-1}U^{\ast }U-Q\right\Vert <\varepsilon )\geq 
1-2\exp \left( -c\varepsilon \min \{n^{1/2}/(\sqrt{D}+\sqrt{D\ln D}),n/(D\ln 
D)\}\right) . 
\end{equation*}%
The probability is not less than $1-\delta $ whenever condition (\ref%
{better_bound}) holds. 
\end{Proof}

Comparing  (\ref{better_bound}) with (\ref{minimum:n}), we have gained on the 
exponent of $D$. However, the quantity $\ln (\delta ^{-1})$ now enters 
quadratically instead of linearly, and an unspecified constant appears in 
the lower bound for $n$.

\subsection{The Non-I.I.D. Case and Other Generalizations\label{noniid}} 

In this section we generalize the results to the case where the random 
matrix $U\in \mathbb{C}^{n\times D}$ has independent rows which, however, 
need not be identically distributed. In the course of this generalization we 
also obtain some slight improvements in the case of i.i.d.~rows discussed 
above. Apart from the assumption of independent rows, we assume that the 
matrix $U$ satisfies the following condition: The moment generating 
functions of the random variables $\Re (\overline{u_{tk}}u_{tj})$ and $\Im (%
\overline{u_{tk}}u_{tj})$ exist for all $1\leq t\leq n$ and $1\leq k,j\leq D$%
; i.e., there exists $x_{0}>0$ such that for all $1\leq t\leq n$ and $1\leq 
k,j\leq D$  
\begin{equation} 
\mathbb{E}\left[ \exp (x\Re (\overline{u_{tk}}u_{tj}))\right] <\infty 
,\qquad \mathbb{E}\left[ \exp (x\Im (\overline{u_{tk}}u_{tj}))\right] <\infty 
\label{moment_gen_func} 
\end{equation}%
holds for all $x<x_{0}$. Note that $x_{0}$ will depend on the distribution 
of $U$ and thus may depend on $n$ and $D$. Furthermore, we set%
\begin{equation*} 
Q^{(t)}:=\mathbb{E(}u_{t\cdot }^{\ast }u_{t\cdot })\in \mathbb{C}^{{D}\times  
{D}} 
\end{equation*}%
with entries $q_{kj}^{(t)}$ and%
\begin{equation} 
Q_{n}\,:=n^{-1}\sum_{t=1}^{n}Q^{(t)}=\mathbb{E}[{n}^{-1}U^{\ast }U]\in  
\mathbb{C}^{{D}\times {D}}.  \label{def:Q} 
\end{equation} 

As in Section~\ref{iid}, assumption (\ref{moment_gen_func}) is seen to be 
equivalent to the existence of finite constants $M_{kj1}^{(t)}\geq 0$, $%
M_{kj2}^{(t)}\geq 0$, $v_{kj1}^{(t)}\geq 0$, $v_{kj2}^{(t)}\geq 0$, such 
that for all $\ell \geq 2$  
\begin{align} 
\mathbb{E}\left[ |\Re (\overline{u_{tk}}u_{tj}-q_{kj}^{(t)})|^{\ell }\right] 
\,& \leq \,2^{-1}\ell !\,(M_{kj1}^{(t)})^{\ell -2}v_{kj1}^{(t)}, 
\label{moment_estimateRe} \\ 
\mathbb{E}\left[ |\Im (\overline{u_{tk}}u_{tj}-q_{kj}^{(t)})|^{\ell }\right] 
\,& \leq \,2^{-1}\ell !\,(M_{kj2}^{(t)})^{\ell -2}v_{kj2}^{(t)} 
\label{moment_estimateIm} 
\end{align}%
hold for all $1\leq t\leq n$ and $1\leq k,j\leq D$. If $%
M_{kji}^{(t)}v_{kji}^{(t)}=0$ then we may assume without loss of generality 
that $M_{kji}^{(t)}=v_{kji}^{(t)}=0$. 

For fixed $n$ it is always possible to choose the constants on the 
right-hand side of (\ref{moment_estimateRe}) and (\ref{moment_estimateIm}) 
independent of $t$. However, for $n\rightarrow \infty $ the resulting 
conditions in Theorem \ref{thm:main_2} below would become unnecessarily 
restrictive in the non-identically distributed case. Furthermore, allowing 
the constants to depend on $k,j$ and to be different in (\ref%
{moment_estimateRe}) and (\ref{moment_estimateIm}), provides some extra 
flexibility which results in an improved, albeit more complex bound even in 
the case of i.i.d.~rows. 

\begin{remark} 
Condition (\ref{moment_estimateRe}) necessarily implies $v_{kj1}^{(t)}\geq 
\sigma _{kj1}^{(t)2}$, where $\sigma _{kj1}^{(t)2}$ denotes the variance of $%
\Re (\overline{u_{tk}}u_{tj}-q_{kj}^{(t)})$. Furthermore, observe that given 
condition (\ref{moment_estimateRe}) is satisfied, it is also always 
satisfied with $v_{kj1}^{(t)}=\sigma _{kj1}^{(t)2}$. [This is obvious if $%
\sigma _{kj1}^{(t)2}=0$, and otherwise follows by replacing $M_{kj1}^{(t)}$ 
with $M_{kj1}^{(t)}v_{kj1}^{(t)}/\sigma _{kj1}^{(t)2}$, observing that $%
v_{kj1}^{(t)}/\sigma _{kj1}^{(t)2}\geq 1$ as noted before.] Similar comments 
apply to condition (\ref{moment_estimateIm}).\hfill $\square $  
\end{remark} 

\begin{remark} 
\label{rem:bounded_2} If the random variables $u_{tk}$ are bounded, i.e.,  
\begin{equation*} 
|\Re (\overline{u_{tk}}u_{tj}-q_{kj}^{(t)})|\,\leq \,C_{kj1}^{(t)}\quad  
\text{and}\quad |\Im (\overline{u_{tk}}u_{tj}-q_{kj}^{(t)})|\,\leq 
\,C_{kj2}^{(t)} 
\end{equation*}%
holds with probability $1$ for all $1\leq t\leq n$, $1\leq k,j\leq D$, then (%
\ref{moment_estimateRe}) and (\ref{moment_estimateIm}) hold with $%
M_{kj1}^{(t)}=C_{kj1}^{(t)}/3$, $M_{kj2}^{(t)}=C_{kj2}^{(t)}/3$, and  
\begin{equation} 
v_{kj1}^{(t)}=\mathbb{E}\left[ (\Re (\overline{u_{tk}}%
u_{tj}-q_{kj}^{(t)}))^{2}\right] \text{ and }v_{kj2}^{(t)}=\mathbb{E}\left[ 
(\Im (\overline{u_{tk}}u_{tj}-q_{kj}^{(t)}))^{2}\right] . 
\label{vij_bounded} 
\end{equation}%
This follows exactly as in Remark \ref{rem:bounded}.\hfill $\square $ 
\end{remark} 

In order to present the generalization of Theorem \ref{thm:main} we introduce%
\begin{equation*} 
v_{kj1n}:=\dsum\limits_{t=1}^{n}v_{kj1}^{(t)},\quad 
v_{kj2n}:=\dsum\limits_{t=1}^{n}v_{kj2}^{(t)} 
\end{equation*}%
and  
\begin{equation*} 
M_{kj1n}=\max \{M_{kj1}^{(t)}:1\leq t\leq n\},\quad M_{kj2n}=\max 
\{M_{kj2}^{(t)}:1\leq t\leq n\}. 
\end{equation*}%
Furthermore, set $v_{n}=\max \{v_{kj1n},v_{kj2n}:1\leq k,j\leq D\}$ and $%
M_{n}=\max \{M_{kj1n},M_{kj2n}:1\leq k,j\leq D\}$. Note that $v_{n}$ and $%
M_{n}$ depend on the distribution of the random matrix $U$ and hence may 
depend on $D$. The expression on the right-hand side of (\ref{prob_bound}) 
below is the direct generalization of (\ref{prob_bound_1}) to the 
non-identically distributed case, whereas the bound $1-\Psi $ given in (\ref%
{union_bound_2}) below is an improvement (even in the case of i.i.d.~rows). 

\begin{Theorem} 
\label{thm:main_2} Assume that the rows $u_{1\cdot },\ldots ,u_{{n\cdot }}$ 
of $U$ are independent random vectors in $\mathbb{C}^{{D}}$ whose entries 
satisfy the moment bounds (\ref{moment_estimateRe}) and (\ref%
{moment_estimateIm}). Then, for every $\varepsilon >0$, the operator norm 
satisfies  
\begin{equation*} 
\left\Vert {n}^{-1}U^{\ast }U-Q_{n}\right\Vert <\varepsilon 
\end{equation*}%
with probability at least $1-\Psi $ where $\Psi $ is defined in (\ref%
{union_bound_2}) below. Furthermore,%
\begin{equation} 
1-\Psi \geq 1-4{D}^{2}\exp \left( -\frac{n\varepsilon ^{2}}{{D}^{2}\left( 
4n^{-1}v_{n}+2\sqrt{2}D^{-1}M_{n}\varepsilon \right) }\right) . 
\label{prob_bound} 
\end{equation}%
In particular, with probability not less than $1-\Psi $, the extremal 
eigenvalues of $n^{-1}U^{\ast }U$ satisfy%
\begin{equation} 
\lambda _{\min }(Q_{n})-\varepsilon <\lambda _{\min }({n}^{-1}U^{\ast 
}U)\leq \lambda _{\max }({n}^{-1}U^{\ast }U)<\lambda _{\max 
}(Q_{n})+\varepsilon \, .  \label{eig_estim} 
\end{equation}%
Consequently the condition number of $U^{\ast }U$ is bounded by $\frac{%
\lambda _{\max }(Q_{n})+\varepsilon }{\lambda _{\min }(Q_{n})-\varepsilon }$ 
with probability not less than $1-\Psi $, provided that $Q_{n}$ defined in (%
\ref{def:Q}) is non-singular and $\varepsilon \in (0,\lambda _{\min 
}(Q_{n})) $. 
\end{Theorem} 

\begin{Proof} 
Exactly as in the proof of Theorem \ref{thm:main} we arrive at  
\begin{align} 
& \mathbb{P}(\left\Vert {n}^{-1}U^{\ast }U-Q_{n}\right\Vert \geq \varepsilon 
)  \label{11} \\ 
& \leq \sum_{k,j=1}^{{D}}\mathbb{P}\left( \left\vert \sum_{t=1}^{{n}}\Re (%
\overline{u_{tk}}u_{tj}-q_{kj}^{(t)})\right\vert \geq \frac{{n}\varepsilon }{%
\sqrt{2}{D}}\right) +\sum_{k,j=1}^{{D}}\mathbb{P}\left( \left\vert 
\sum_{t=1}^{{n}}\Im (\overline{u_{tk}}u_{tj}-q_{kj}^{(t)})\right\vert \geq  
\frac{{n}\varepsilon }{\sqrt{2}{D}}\right) .  \notag 
\end{align}%
Again using inequality (\ref{bennett}) for each $k,j$ gives  
\begin{equation*} 
\mathbb{P}\left( \left\vert \sum_{t=1}^{{n}}\Re (\overline{u_{tk}}%
u_{tj}-q_{kj}^{(t)})\right\vert \geq \frac{{n}\varepsilon }{\sqrt{2}{D}}%
\right) \leq 2\exp \left( -\frac{n\varepsilon ^{2}}{{D}^{2}\left( 
4n^{-1}v_{kj1n}+2\sqrt{2}D^{-1}M_{kj1n}\varepsilon \right) }\right) 
\end{equation*}%
and similarly for the imaginary part. Hence, we finally obtain $\mathbb{P}%
\left( \Vert {n}^{-1}U^{\ast }U-Q\Vert \geq \varepsilon \right) \,\leq \Psi $ 
where%
\begin{eqnarray} 
\Psi &=&2\sum_{i=1}^{2}\sum_{k,j=1}^{{D}}\exp \left( -\frac{n\varepsilon ^{2}%
}{{D}^{2}\left( 4n^{-1}v_{kjin}+2\sqrt{2}D^{-1}M_{kjin}\varepsilon \right) }%
\right)  \label{union_bound_2} \\ 
&\leq &4{D}^{2}\exp \left( -\frac{n\varepsilon ^{2}}{{D}^{2}\left( 
4n^{-1}v_{n}+2\sqrt{2}D^{-1}M_{n}\varepsilon \right) }\right) .  \notag 
\end{eqnarray} 
\end{Proof} 

\begin{remark} 
A sufficient condition for $Q_{n}$ to be non-singular is that at least one 
of the matrices $Q^{(t)}$ has this property. The argument in Remark \ref%
{nonsing} shows that the latter is the case if the distribution of $u_{t.}$ 
is not concentrated on a $(D-1)$-dimensional linear subspace of $\mathbb{C}%
^{D}$. However, note the possibility that nevertheless $\lambda _{\min 
}(Q_{n})\rightarrow 0$ as $n\rightarrow \infty $.\hfill $\square $ 
\end{remark} 

\begin{remark} 
\label{zero}(i) In case the $(k,j)$-element of ${n}^{-1}U^{\ast }U-Q_{n}$ is 
zero with probability $1$, the corresponding terms on the right-hand side of 
(\ref{11}) are zero and do not contribute to the bound in (\ref{11}). Due to 
the independence assumption, the $(k,j)$-element is zero if and only if $%
\overline{u_{tk}}u_{tj}-q_{kj}^{(t)}=0$ with probability $1$ for every $t$. 
Hence, we may set $v_{kj1}^{(t)}=v_{kj2}^{(t)}=M_{kj1}^{(t)}=M_{kj2}^{(t)}=0$ 
which shows that the corresponding terms in the bound $\Psi $ are also 
automatically zero. However, in this case the bound (\ref{11}) and the 
subsequent bounds can be improved in that in the $(k,j)$-th term in both 
sums on the right-hand side of (\ref{11}) the constant $D$ can be replaced 
by $D_{k}$, where $D_{k}$ denotes the number of non-zero elements in the $k$%
-th row of ${n}^{-1}U^{\ast }U-Q_{n}$. 

(ii) A similar remark applies in the case that some or all elements of ${n}%
^{-1}U^{\ast }U-Q_{n}$ are real (or imaginary). Cf. Remark \ref{real}. 

\end{remark} 

\section{Random Sampling of Trigonometric Polynomials Revisited\label{Sec_4}} 

We now return to the special case of sampling trigonometric polynomials on 
uniformly distributed random points and show how the results in the previous 
sections can be improved. The analysis is based on techniques developed in  
\cite{ra05-7} for the recovery of sparse trigonometric polynomials from 
random samples by basis pursuit ($\ell _{1}$-minimization) and orthogonal 
matching pursuit. Some of the ideas are inspired by the pioneering work of 
Cand{\`{e}}s, Romberg and Tao in \cite{carota06}. 

\begin{Theorem} 
\label{thm_uniform} Let $\Gamma \subset \mathbb{Z}^{d}$ of size $|\Gamma |={D%
}$ and let $x_{1},\ldots ,x_{{n}}$ be i.i.d.~random variables that are 
uniformly distributed on $[0,1]^{d}$. Let $U$ be the associated random 
Fourier matrix given by (\ref{defU}). Choose $0<\varepsilon <1$, $0<\alpha 
<\varepsilon ^{2}$, and $\delta >0$. If  
\begin{equation} 
\left\lfloor \frac{\alpha {n}}{3{D}}\right\rfloor \,\geq \,\left[ \ln \left(  
\frac{\varepsilon ^{2}}{\alpha }\right) \right] ^{-1}\ln \left( \frac{{D}}{%
\delta (1-\alpha )}\right) ,  \label{cond_nD} 
\end{equation}%
then, with probability at least $1-\delta $, we have  
\begin{equation*} 
\left\Vert {n}^{-1}U^{\ast }U-I\right\Vert <\varepsilon 
\end{equation*}%
and hence%
\begin{equation*} 
1-\varepsilon <\lambda _{\min }({n}^{-1}U^{\ast }U)\leq \lambda _{\max }({n}%
^{-1}U^{\ast }U)<1+\varepsilon . 
\end{equation*}%
Consequently, the condition number of $U^{\ast }U$ is bounded by $\frac{%
1+\varepsilon }{1-\varepsilon }$ with probability $\geq 1-\delta $. 
\end{Theorem} 

For instance, the choice $\alpha =\varepsilon ^{2}/e$ gives  
\begin{equation} 
n\geq \frac{3De}{\varepsilon ^{2}}\left[ \ln \left( \frac{D}{\delta }\right) 
+2-\ln (e-1)\right]  \label{minimum:n2} 
\end{equation}%
as a simple sufficient condition. 

Comparing (\ref{minimum:n2}) with (\ref{minimum:n}) or (\ref{minimum:nb}), we 
have gained on the exponent in $D$; compared with Theorem \ref{thm:detprob} 
and (\ref{est:detprob}), the constants are now independent of the dimension $%
d$ of the state space ($=$ the number of variables); compared with (\ref%
{better_bound}), the term $\ln (\delta ^{-1})$ only enters linearly in (\ref%
{cond_nD}) and (\ref{minimum:n2}) instead of quadratically and there is now 
no restriction on $\Gamma $. Moreover, the constants are explicit and small. 

\subsection{Proof of Theorem \protect\ref{thm_uniform}} 

We introduce the polynomials  
\begin{equation} 
F_{m}(z)\,=\,\sum_{k=1}^{\lfloor m/2\rfloor }S_{2}(m,k)z^{k},\qquad m\in  
\mathbb{N},  \label{def_Fn0} 
\end{equation}%
where $S_{2}(m,k)$ are the associated Stirling numbers of the second kind. 
These are connected to the combinatorics of certain set partitions, and they 
can be computed by means of their exponential generating function, see \cite[%
formula (27), p.77]{ri58} or Sloane's A008299 in \cite{slXX},  
\begin{equation} 
\sum_{m=1}^{\infty }F_{m}(z)\frac{x^{m}}{m!}\,=\,\exp \big(z(e^{x}-x-1)\big). 
\label{def_Fn} 
\end{equation}%
Further, we define  
\begin{equation} 
G_{m}(z)\,:=\,z^{-m}F_{m}(z).  \label{def_Gn} 
\end{equation}%
Using the $G_m$'s, we first establish a more general result from which 
Theorem \ref{thm_uniform} will follow. 

\begin{Theorem} 
\label{thm_Gm}Let $\Gamma \subset \mathbb{Z}^{d}$ of size $|\Gamma |={D}$ 
and let $x_{1},\ldots ,x_{{n}}$ be i.i.d.~random variables that are 
uniformly distributed on $[0,1]^{d}$. Let $U$ be the associated random 
Fourier matrix given by (\ref{defU}), and let $\varepsilon >0$. Then, for 
every $m\in \mathbb{N}$, we have%
\begin{equation*} 
\left\Vert {n}^{-1}U^{\ast }U-I\right\Vert <\varepsilon , 
\end{equation*}%
and hence  
\begin{equation*} 
1-\varepsilon <\lambda _{\min }({n}^{-1}U^{\ast }U)\,\leq \,\lambda _{\max }(%
{n}^{-1}U^{\ast }U)\,<\,1+\varepsilon \, , 
\end{equation*}%
with probability at least%
\begin{equation*} 
1-\varepsilon ^{-2m}{D}G_{2m}({n}/{D}). 
\end{equation*} 
\end{Theorem} 

\begin{Proof} 
Again, the estimates for the eigenvalues follow from the inequality $\Vert {n%
}^{-1}U^{\ast }U-I\Vert <\varepsilon $. Furthermore, since ${n}^{-1}U^{\ast 
}U-I$ is self-adjoint, we have for every $m\in \mathbb{N}$  
\begin{equation*} 
\Vert {n}^{-1}U^{\ast }U-I\Vert \,=\,\Vert ({n}^{-1}U^{\ast }U-I)^{m}\Vert 
^{1/m}\,\leq \,\Vert ({n}^{-1}U^{\ast }U-I)^{m}\Vert _{F}^{1/m}, 
\end{equation*}%
where $\Vert \cdot \Vert _{F}$ denotes the Frobenius norm, $\Vert A\Vert 
_{F}\,=\,\sqrt{\func{Tr}(AA^{\ast })}$. Consequently,%
\begin{equation*} 
\mathbb{P}(\Vert {n}^{-1}U^{\ast }U-I\Vert \geq \varepsilon )\leq \mathbb{P}%
(\Vert ({n}^{-1}U^{\ast }U-I)^{m}\Vert _{F}\geq \varepsilon ^{m}).\, 
\end{equation*}%
We now apply Markov's inequality and obtain that  
\begin{equation*} 
\mathbb{P}(\Vert ({n}^{-1}U^{\ast }U-I)^{m}\Vert _{F}\geq \varepsilon 
^{m})\leq \varepsilon ^{-2m}\mathbb{E}\left[ \Vert ({n}^{-1}U^{\ast 
}U-I)^{m}\Vert _{F}^{2}\right] . 
\end{equation*}%
The latter expectation was studied in \cite[Section 3.3]{ra05-7}, see also 
Lemma 3.3 in \cite{ra05-7}: It was shown that  
\begin{equation} 
\mathbb{E}\left[ \Vert ({n}^{-1}U^{\ast }U-I)^{m}\Vert _{F}^{2}\right] 
\,\leq \,{D}\,G_{2m}({n}/{D}),  \label{H0_Frob} 
\end{equation}%
which concludes the proof. 
\end{Proof} 

We now show how Theorem \ref{thm_uniform} follows from Theorem~\ref{thm_Gm}. 
This is done by estimating $G_m$ and by a diligent choice of the free 
parameter $m$. 

We set the oversampling rate to be $\theta =n/{D}$. In \cite[Section 3.5]%
{ra05-7} it was shown that  
\begin{equation*} 
G_{2m}(\theta )\,\leq \,(3m/\theta )^{m}\frac{1-(3m/\theta )^{m}}{%
1-(3m/\theta )}. 
\end{equation*}%
For given $0<\alpha <1$ and $\theta $, we choose $m=m(\theta )\in \mathbb{N}$ 
such that $(3m(\theta )/\theta )\leq \alpha <1$. Note that this is possible 
since $\left\lfloor \alpha n/3D\right\rfloor \geq 1$ follows from the 
assumptions of the theorem. In the following we will take the value  
\begin{equation*} 
m(\theta )\,=\,\left\lfloor \alpha \theta /3\right\rfloor \,, 
\end{equation*}%
and obtain that  
\begin{equation} 
G_{2m(\theta )}(\theta )\,\leq \,\frac{\alpha ^{m(\theta )}}{1-\alpha }. 
\label{eq:c5} 
\end{equation}%
In view of Theorem~\ref{thm_Gm} we want to achieve $\varepsilon ^{-2m}{D}%
\,G_{2m}(\theta )\leq \delta $. By (\ref{eq:c5}) this inequality is 
satisfied if  
\begin{equation*} 
{D}\varepsilon ^{-2m(\theta )}\frac{\alpha ^{m(\theta )}}{1-\alpha }\,\leq 
\,\delta , 
\end{equation*}%
which is equivalent to  
\begin{equation*} 
\ln \Big(\frac{\varepsilon ^{2}}{\alpha } \Big)\, m(\theta )\,\geq \,\ln %
\Big( \frac{D}{(1-\alpha ) \delta}\Big) . 
\end{equation*}%
Since $\alpha <\varepsilon ^{2}$ by assumption, Theorem \ref{thm_uniform} 
follows. 

Let us mention that an  estimate of the form $n \geq C D \log D \epsilon ^{-2}
\delta ^{-2}$ can also be derived from Rudelson's estimates for random
vectors in isotropic position~\cite{Ru99}. 

\section*{Acknowledgments} 

We would like to thank Richard Nickl for helpful discussions.


\end{document}